\documentstyle[12pt]{article}

\date{ }
\begin{document}
\setcounter{page}{1}
\rightline{q-alg/9501021}
\rightline{CPTH-A344.0195}
\begin{center}
\begin{Large}
The fundamental invariant of the Hecke algebra $H_n(q)$ characterizes the
representations of $H_n(q)$, $S_n$, $SU_q(N)$ and $SU(N)$
\vspace{10pt}

$\mbox{J. Katriel}^{*}$
, B. Abdesselam and A. Chakrabarti \\
\vspace{10pt}
{\small \sl Centre de Physique Th\'eorique}\\
{\small \sl Ecole Polytechnique}\\
{\small \sl 91128 Palaiseau Cedex, France}
\end{Large}
\end{center}
\vspace*{105mm}
\ \hrulefill\ \hfill \

\begin{small}
${\ }^{*}$Permanent address: Department of Chemistry,
Technion, 32000 Haifa, Israel.
\end{small}

\newpage
\pagestyle{plain}

\vspace*{2 cm}

\begin{abstract}

The irreducible representations (irreps) of the Hecke algebra $H_n(q)$
are shown to be completely
characterized by the fundamental invariant of this algebra, $C_n$.
This fundamental invariant is related to the quadratic Casimir operator,
${\cal{C}}_2$, of $SU_q(N)$, and reduces to the transposition class-sum,
$[(2)]_n$, of $S_n$ when $q\rightarrow 1$.
The projection operators constructed in terms of $C_n$ for the
various irreps of $H_n(q)$ are well-behaved in the limit $q\rightarrow 1$,
even when approaching degenerate eigenvalues of $[(2)]_n$.
In the latter case, for which the
irreps of $S_n$ are not fully characterized by the
corresponding eigenvalue of the transposition class-sum, the limiting form of
the projection
operator constructed in terms of
$C_n$ gives rise to factors that depend on
higher class-sums of $S_n$, which effect the desired characterization.
Expanding this limiting form of the projection operator into a linear
combination of class-sums of $S_n$, the coefficients constitute the
corresponding row in the character table of $S_n$.
The properties of the fundamental invariant are used to formulate a simple
and efficient recursive procedure for the evaluation of the traces of
the Hecke algebra. The closely related quadratic Casimir operator of $SU_q(N)$
plays a similar role, providing a complete characterization of the irreps
of $SU_q(N)$ and - by constructing appropriate projection operators and then
taking the $q\rightarrow 1$ limit - those of $SU(N)$ as well, even when the
quadratic
Casimir operator of the latter does not suffice to specify its irreps.

\end{abstract}

\newpage

\section{Introduction}

\hspace*{6 mm}
The ordinary representation theory of the symmetric group, whose centenial
is soon to be celebrated, is not the most likely arena for the emergence
of hitherto unnoticed simplifications.
The interest in the Hecke algebra $H_n(q)$, that reduces to the group algebra
of the symmetric group in the limit $q\rightarrow 1$, is much more recent.
It has been pursued as a vehicle for the construction of representations of
the braid group and also in the context of the recent study of the quantum
unitary groups, with
respect to which $H_n(q)$ has been shown to play the role of the q-analogue
of the Weyl group.
For recent expositions of the mathematical structure and of the physical
relevance of the Hecke algebra, jointly providing access to many earlier
references,
see refs. [1-11].

The central theme of the present article is that the Hecke algebra point of
view
provides a remarkable conceptual simplification in the formulation of the
representation theory of the symmetric group.
More specifically, the fundamental invariant of $H_n(q)$, that in the limit
$q\rightarrow 1$ reduces to the transposition class-sum of $S_n$, is shown
to be sufficient to effect all the applications that conventionally require
the use of the full character table of $S_n$. In fact, the character
table of $H_n(q)$, and, {\it a fortiori,} that of $S_n$, are shown to be
extractable
by appropriate manipulations of that fundamental invariant.

A completely analogous situation holds true for the
fundamental invariant of $SU_q(N)$, its quadratic Casimir operator.
It turns out that the irreducible representations (irreps) of $SU_q(N)$,
and, {\it a fortiori,}
of $SU(N)$, are fully characterized by the quadratic Casimir of the former.
In fact, the eigenvalues of the quadratic Casimir operator of $SU_q(N)$
can be expressed in terms of those of the fundamental invariant of $H_n(q)$ in
a manner that is the q-analogue of the relation noted by Gross \cite{Gross}
between the eigenvalues of the quadratic Casimir operator of $SU(N)$ and the
transposition class-sum of $S_n$.

The structure of the present article is as follows:
In section 2 we review some pertinent features of the representation theory of
the
symmetric group. In section 3 we show that the fundamental invariant of the
Hecke algebra
fully characterizes its irreps.
In section 4 we construct projection operators for the irreps of $H_n(q)$ in
terms of its
fundamental invariant, and study their limiting forms for $q\rightarrow 1$,
that
yield the familiar group-theoretical projection operators for the irreps of
$S_n$.
In section 5 we discuss the extraction of the traces of the Hecke algebra from
its fundamental invariant, and present an efficient  procedure for the
evaluation of
these traces.
The relation between the quadratic Casimir operator of $SU_q(N)$ and the
fundamental
invariant of $H_n(q)$ is presented in section 6, where, for generic $q$,
the former is shown to characterize all the corresponding irreps.
Some concluding remarks are made in section 7.

\section{Class-sums and classification of irreps of $S_n$}

\hspace*{6 mm}
Both the irreps and the conjugacy classes
of the symmetric group
$S_n$ can be labeled by means of partitions of $n$.
We shall denote irreps by bracketed sets of bold face
integers specifying
the row lengths of the corresponding Young diagram,
and conjugacy classes by sets of individually parenthesized integers
that specify the cycle-structure characterising the class.
A class-sum will be denoted by bracketing the class-symbol, appending a suffix
that
specifies $n$, and usually suppressing cycles of unit length.

As is well-known, the irreps are fully characterized by means of the
set of eigenvalues of the class-sums.
In fact, it has been shown by Kramer \cite{Kramer} that the single-cycle
class-sums,
$\{ [(p)]_n \, ;\, p=2,\, 3,\, \cdots,\, n\}$, generate the center of the
corresponding group-algebra and are therefore sufficient to characterize the
irreps.
Moreover, it has been noted that the subset of $k$ class-sums $\{ [(2)]_n,\,
[(3)]_n,\,
\cdots,\, [(k+1)]_n\}$ is sufficient to characterize the irreps of all
symmetric groups
$S_n$ with $n$ not larger than some $n_{\max}(k)$, that is considerably larger
than $k$.
Thus, $n_{\max}(1)=5$,
the first instance in which the transposition class-sum has the same eigenvalue
for
two different irreps being the irreps
${\bf{[4,1,1]}}$ and ${\bf{[3,3]}}$
of $S_6$. Similarly,
$n_{\max}(2)=14$, $n_{\max}(3)=23$ and $n_{\max}(4)=41$. The
rather conservative bound $n_{\max}(k)<2^{2k+1}$ has been found to hold
\cite{Katriel95}
for arbitrary $k$.

The eigenvalues of the class-sums of the symmetric group,
{\it i.e.,} its central characters, have been shown
in ref. \cite{Katriel1} to be
expressible as polynomials in the symmetric power-sums
$$\sigma_k^{\Gamma}=\sum_{(i,j)\in{\Gamma}} (j-i)^k \; .$$
Here, $\Gamma$ stands for a Young diagram and $(i,j)$ are the row and column
indices
of the boxes in this diagram, counted matrixlike. The difference $j-i$ has been
referred to as the content of the box $(i,j)$. Algorithms for the
construction of expressions for the central characters in terms of these
symmetric power-sums were
presented in refs. \cite{Katriel2} and \cite{Katriel3} for single-cycle and
arbitrary class-sums, respectively.

For the first four single-cycle class-sums the eigenvalues are given by
\begin{eqnarray}
\label{equation:ei}
\lambda_{[(2)]_n}^{\Gamma}&=&\sigma_1^{\Gamma} \\
\lambda_{[(3)]_n}^{\Gamma}&=&\sigma_2^{\Gamma}-\frac{n(n-1)}{2} \nonumber \\
\lambda_{[(4)]_n}^{\Gamma}&=&\sigma_3^{\Gamma}-(2n-3)\sigma_1^{\Gamma}
\nonumber \\
\lambda_{[(5)]_n}^{\Gamma}&=&\sigma_4^{\Gamma}-(3n-10)\sigma_2^{\Gamma}
-2(\sigma_1^{\Gamma})^2+\frac{n(n-1)(5n-19)}{6} \nonumber
\end{eqnarray}

\noindent
The expressions for $\lambda_{[(2)]_n}^{\Gamma}$ and
$\lambda_{[(3)]_n}^{\Gamma}$ were originally given by Jucys \cite{Jucys} and by
Suzuki \cite{Suzuki}.

\section{The Hecke algebra $H_n(q)$}

\hspace*{6 mm}
The Hecke algebra $H_n(q)$ is defined in terms of the generators
$g_1,\, g_2, \, \cdots, \, g_{n-1}$ and the relations
\begin{equation}
\label{equation:Hecke}
\begin{array}{ll}
g_i^2=(q-1)g_i+q &\;\;\; i=1,\, 2,\, \cdots,\, n-1 \\
g_ig_{i+1}g_i=g_{i+1}g_ig_{i+1} &\;\;\; i=1,\, 2,\, \cdots,\, n-2 \\
g_ig_j=g_jg_i &\;\;\; {\mbox{if }} |i-j|\geq 2
\end{array}
\end{equation}

For $q=1$ the relations specified above reduce to the generating relations of
the
symmetric group, $S_n$. In particular, $g_i$ reduces to the transposition
$(i,i+1)$.

When $q$ is neither zero nor a $k$th root of unity, $k=2,\, 3,\, \cdots, \, n$,
the irreps of $H_n(q)$ are labelled by Young diagrams with $n$ boxes
\cite{Wenzl,Vershik}.

The fundamental invariant is
\begin{eqnarray}
\label{equation:fund}
C_n &=& g_1+g_2+\cdots +g_{n-1}+
\frac{1}{q}(g_1g_2g_1+g_2g_3g_2+\cdots +g_{n-2}g_{n-1}g_{n-2}) \nonumber \\ &+&
\frac{1}{q^2}(g_1g_2g_3g_2g_1+g_2g_3g_4g_3g_2+\cdots
+g_{n-3}g_{n-2}g_{n-1}g_{n-2}g_{n-3})
\nonumber \\ &+&\cdots \nonumber \\ &+&
\frac{1}{q^{n-2}}g_1g_2\cdots g_{n-2}g_{n-1}g_{n-2}\cdots g_2g_1
\end{eqnarray}

$C_n$ is, up to an overall factor of $q$, the sum of the Murphy operators
introduced
by Dipper and James \cite{Dipper}. Using their corollary 2.3 it follows that
$C_n$
belongs to the center of $H_n(q)$, {\it i.e.,} commutes with all the generators
$g_i, \;\; i=1,\, 2,\, \cdots,\, n-1$.
As an illustration consider $C_3=g_1 + g_2 + \frac{1}{q}g_1g_2g_1$.
In this case \hfill\break
$[g_1,C_3]=g_1g_2-g_2g_1 +\frac{1}{q}\Big((q-1)g_1+q\Big)g_2g_1
-\frac{1}{q}g_1g_2\Big( (q-1)g_1+q\Big)=0$.

Using the relation between $C_n$ and the Murphy operators it follows from
corollary 3.13
in ref. \cite{Dipper} that
the eigenvalues of the fundamental invariant are given by the expression
\begin{equation}
\label{equation:eigen}
\Lambda_n^{\Gamma}(q)=q\sum_{(i,j)\in\Gamma} [j-i]_q
\end{equation}

\noindent
where $[k]_q=\frac{q^k-1}{q-1}$.
We shall refer to $q[j-i]_q$ as the q-content of the box $(i,j)$.

For $q\rightarrow 1$ the fundamental invariant reduces to the symmetric-group
transposition class-sum,
$[(2)]_n$, and eq. \ref{equation:eigen} reduces to the expression for the
corresponding eigenvalues, \hfill\break
$\lambda_{[(2)]_n}^{\Gamma}=\sigma_1^{\Gamma}$ ({\it cf.} eq.
\ref{equation:ei}).

To illustrate the expression for the eigenvalues consider $C_2=g_1$.
The corresponding eigenvalues are the roots of the Cayley equation
$\Lambda^2=(q-1)\Lambda+q$ that is obtained from $g_1^2=(q-1)g_1+q$ by
replacing
$g_1$ by its eigenvalue $\Lambda$. The roots of this quadratic equation,
$\Lambda=q,\, -1$, coincide with $\Lambda_{C_2}^{[2]}(q)$ and
$\Lambda_{C_2}^{[1,1]}$,
respectively, evaluated using eq. \ref{equation:eigen}.
A less trivial illustration is provided by $H_3(q)$. Here,
$$C_3=g_1+g_2+\frac{1}{q}g_1g_2g_1\; ,$$

\noindent
so that

$$C_3^2=3q+2(q-1)C_3+(q+1+\frac{1}{q})D_3$$

\noindent
where
$D_3=g_1g_2+g_2g_1+\frac{q-1}{q}g_1g_2g_1 \; ,$

\noindent
and

$$C_3^3=(q^3-1)+q(q^2+q+5+\frac{1}{q}+\frac{1}{q^2})C_3
+2(q-1)C_3^2 +(q^3+2q^2+q-1-\frac{2}{q}-\frac{1}{q^2})D_3$$

\noindent
Eliminating $D_3$ and replacing $C_3$ by its eigenvalue we obtain the Cayley
equation
$$\Lambda^3-(q-1)(q+4+\frac{1}{q})\Lambda^2
+(q^3-q^2-9q-1+\frac{1}{q})\Lambda+(q-1)(2q^2+5q+2)=0$$

\noindent
whose roots, $\Lambda=q^2+2q,\; q-1,\; -2-\frac{1}{q}$, agree with the values
obtained
for $\Lambda_{C_3}^{[3]},\; \Lambda_{C_3}^{[2,1]},$ and
$\Lambda_{C_3}^{[1,1,1]}$,
respectively, using eq. \ref{equation:eigen}.

The most important property of the eigenvalues of $C_n$ from the point of view
of the
applications to be
presented below is that, being polynomials in $q$ and $\frac{1}{q}$,
for generic $q$ they obtain distinct
values for the various irreps of $H_n(q)$.
As an illustration consider the two irreps ${\bf [4,1,1]}$ and ${\bf [3,3]}$ of
$H_6(q)$,
whose $S_6$ counterparts both correspond to the eigenvalue $3$ of $[(2)]_6$.
Using eq. \ref{equation:eigen} we obtain
$\Lambda_{C_6}^{[4,1,1]}=q^3+2q^2+3q-2-\frac{1}{q}$
and
$\Lambda_{C_6}^{[3,3]}=q^2+3q-1$.
These two eigenvalues are equal to one another only when $q$ is either a
sqaure- or
a cubic root of unity. Hence, while they are equal when $q=1$, there is some
neighbourhood of $q=1$ within which they are distinct.

To see this property in the general setting we note that
the principal diagonal of a Young diagram consists of boxes with equal row and
column indices, {\it i.e.,} of the boxes whose contents are zero.
The boxes with a common content $k$ lie along a diagonal that is parallel to
the
principal diagonal at a distance $k$ above or below it, depending on the sign
of $k$.
Denoting the number of boxes with content $k$ in the Young diagram $\Gamma$ by
$\beta_k^{\Gamma}$, the eigenvalue of $C_n$ can be written in the form
\begin{eqnarray}
\label{equation:boxes}
\Lambda_{C_n}^{\Gamma}(q) &=& \sum_{k>0} (q+q^2+\cdots +q^k)\beta_k^{\Gamma}
-\sum_{k<0} (1+\frac{1}{q}+\cdots +\frac{1}{q^{-k-1}})\beta_k^{\Gamma}
\nonumber \\ &=&
\sum_{k>0}q^k\sum_{\ell \geq k} \beta_{\ell}^{\Gamma}
- \sum_{k<0} \frac{1}{q^{-k-1}} \sum_{\ell \leq k}\beta_{\ell}^{\Gamma}
\nonumber \\ &=&
\sum_{k>0}q^k\pi_k^{\Gamma}
-\sum_{k<0} \frac{1}{q^{-k-1}}\nu_k^{\Gamma}
\end{eqnarray}
{\it i.e.,} $\pi_k^{\Gamma}$, the coefficient of $q^k$, $k>0$,
in $\Lambda_{C_n}^{\Gamma}(q)$, is equal
to the number of boxes with contents larger than or equal to $k$,
and $\nu_k^{\Gamma}$, the coefficient of $\frac{1}{q^{-k-1}}$, $k<0$,
is equal to the number of boxes with contents smaller than or equal to $k$.

Denoting the largest and smallest contents of a given Young diagram by $k_M$
and $k_m$,
respectively, we note that the former must correspond to a single box, located
at the
right end of the top row, and the latter to a single box located at the bottom
of the first column. Thus, $\beta_{k_M}^{\Gamma}=\beta_{k_m}^{\Gamma}=1$.
For $0<k<k_M$ we have from eq. \ref{equation:boxes}
$$\beta_k^{\Gamma}=\pi_k^{\Gamma}-\sum_{\ell>k}\beta_{\ell}^{\Gamma}
=\pi_k^{\Gamma}-\pi_{k+1}^{\Gamma}$$
and for $k_m<k<0$
$$\beta_k^{\Gamma}=\nu_k^{\Gamma}-\sum_{\ell<k}\beta_{\ell}^{\Gamma}=
\nu_k^{\Gamma}-\nu_{k-1}^{\Gamma}$$
Finally, the number of boxes along the principal diagonal (boxes whose contents
are zero)
is $n-\Big(\sum_{k>0}\beta_k^{\Gamma}+\sum_{k<0}\beta_k^{\Gamma}\Big)=
n-(\pi_1^{\Gamma}+\nu_{-1}^{\Gamma})$.

In conclusion, given an eigenvalue of $C_n$ as a polynomial in $q$ and
$\frac{1}{q}$,
the structure of the corresponding Young diagram, expressed in terms of the
lengths
of its diagonals, can readily be recovered.
This conclusion presupposes that different integral powers of $q$ with positive
powers
up to and including $k_M$ and negative powers up to an including $k_m$ are
distinct.
For a Young diagram with $n$ boxes $k_M \leq n-1$ and $k_m \geq -(n-1)$.
Therefore, to assure that different relevant powers of $q$ are distinct, it is
sufficient to require that $q$ is not a root of unity of order less than or
equal to $n$.

It will be convenient to define the operator ${\tilde{C}}_n \equiv
\frac{q-1}{q}C_n$
whose eigenvalues are \hfill\break
$\Lambda_{{\tilde{C}}_n}^{\Gamma}(q)=\sum_{(i,j)\in \Gamma} (q^{j-i}-1)$.

Substituting $q=\exp(\delta)$ we obtain
\begin{equation}
\label{equation:tild}
\Lambda_{{\tilde{C}}_n}^{\Gamma}(\exp(\delta))
= \sum_{k=1}^{\infty} \frac{\delta^k}{k!}\sigma_k^{\Gamma}
\end{equation}
Hence, using eq. \ref{equation:ei},
\begin{equation}
\label{equation:expan}
\Lambda_{{\tilde{C}}_n}^{\Gamma}(\exp(\delta)) =
\delta\lambda_{[(2)]_n}^{\Gamma}
+\frac{\delta^2}{2!}\Big(\lambda_{[(3)]_n}^{\Gamma}+\frac{n(n-1)}{2}\Big)
+\frac{\delta^3}{3!}\Big(\lambda_{[(4)]_n}^{\Gamma}+(2n-3)\lambda_{[(2)]_n}^{\Gamma}\Big)
+\cdots
\end{equation}

In fact, eq. \ref{equation:tild} implies that, for generic $q$, the eigenvalue
$\Lambda_{C_n}^{\Gamma}(q)$
determines all the symmetric power sums
$\{ \sigma_k^{\Gamma}\, | \, k=1,\, 2,\, \cdots, n-1\}$
and, consequently, the central characters of all the class-sums of the
corresponding
symmetric group.

Eq. \ref{equation:expan} implies the correspondence
\begin{eqnarray}
\label{equation:corres}
& & {\tilde{C}}_n\rightarrow
\delta [(2)]_n+\frac{\delta^2}{2!}\Big([(3)]_n+\frac{n(n-1)}{2}\Big)
\\ & & {\phantom{C_2\rightarrow [(2)]_n+\frac{\delta}{2}\Big([(3)]_n}}
+\frac{\delta^3}{3!}\Big([(4)]_n+(2n-3)[(2)]_n
\Big)+\cdots \nonumber
\end{eqnarray}
that should be understood to imply that the operator on the right has the same
eigenvalue for a Young diagram corresponding to an irrep of $S_n$ that
${\tilde{C}}_n$ has for the
same Young diagram taken as an irrep of $H_n(q)$.
This correspondence will be used in the next section.

\section{Projection operators}

\hspace*{6 mm}
The group-theoretical projection operators for the subspace carrying any
desired
irrep can easily be written down in terms of the characters corresponding to
that
representation.
This procedure appears to suggest that full knowledge of the character table is
required,
or, more specifically, the projection operator corresponding to any particular
irrep requires the use of the complete row of characters belonging to that
irrep.

In fact, when a subset of class-sums is sufficient to characterize the irreps,
one only needs the eigenvalues corresponding to these class-sums.
As a very simple illustration consider $S_3$, that has three irreps,
corresponding to
the Young diagrams ${\bf [3]}$, ${\bf [2,1]}$, and ${\bf [1,1,1]}$. The
eigenvalues of
the transposition class-sum corresponding to these irreps are easily evaluated
using
the first of eqs. \ref{equation:ei}, obtaining the respective values 3, 0, and
-3.
The projection operator for any of the three irreps can now be constructed in
terms
of the transposition class-sum. Thus,
$$P_{[2,1]}=\frac{([(2)]_3-3)}{(0-3)}\cdot \frac{([(2)]_3+3)}{(0+3)}$$
where the left factor annihilates the irrep ${\bf [3]}$ and the right factor
annihilates ${\bf [1,1,1]}$. Expanding the product and using the identity
$[(2)]_3[(2)]_3=3+3[(3)]_3$ we obtain
$$P_{[2,1]}=\frac{1}{3}\left( 2-[(3)]_3 \right) \; .$$
By comparison with the group-theoretical form of the projection operator this
expression
provides the row in the character table corresponding to the irrep ${\bf
[2,1]}$,
{\it i.e.,} $\chi((1)^3)=2$, $\chi((1)(2))=0$, and $\chi((3))=-1$.

Let us denote the set of irreps of $S_n$ (and eventually of $H_n(q)$) by
${\cal{IRR}}_n$.
For any irrep $\Gamma_0$ that satisfies
\begin{equation}
\label{equation:cond1}
\lambda_{[(2)]_n}^{\Gamma_0} \neq \lambda_{[(2)]_n}^{\Gamma}
\end{equation}
for all
$\Gamma\in {\cal{IRR}}_n\setminus\{ \Gamma_0 \}$
we write the projection opeartor in the form
$$P_{\Gamma_0}=\prod_{\Gamma\in{\cal{IRR}}_n\setminus\{\Gamma_0\} }
\frac{[(2)]_n-\lambda_{[(2)]_n}^{\Gamma}}
{\lambda_{[(2)]_n}^{\Gamma_0}-\lambda_{[(2)]_n}^{\Gamma} }$$

The procedure illustrated above needs a slight modification for irreps that
do not satisfy eq. \ref{equation:cond1}, {\it i.e.,} are
not fully characterized by the corresponding eigenvalue of the transposition
class-sum. Thus, the two irreps ${\bf [4,1,1]}$ and ${\bf [3,3]}$
of $S_6$ both correspond to the eigenvalue 3 for $[(2)]_6$.
These two irreps can be distinguished
by means of the class-sum $[(3)]_6$, for which their eigenvalues are
4 and -8, respectively.
Furthermore, no other irrep of $S_6$ corresponds to the same eigenvalue of
$[(2)]_6$.
Therefore, to construct a projection operator for any of the above two irreps
one can
first annihilate all the irreps but these two by means of an operator,
$P^{\ast}$, that depends
only on the transposition class-sum.
Explicitly,
$$P^{\ast}=\prod_{\Gamma\in{\cal{IRR}}_6^{\ast}}
\frac{[(2)]_6-\lambda_{[(2)]_6}^{\Gamma}}
{\lambda_{[(2)]_6}^{\ast}-\lambda_{[(2)]_6}^{\Gamma} }$$
where
$\lambda_{[(2)]_6}^{\ast}=\lambda_{[(2)]_6}^{[4,1,1]}=\lambda_{[(2)]_6}^{[3,3]}=3$
and ${\cal{IRR}}_6^{\ast}={\cal{IRR}}_6\setminus\{
{\mbox{\bf [4,1,1]}},\, {\mbox{\bf [3,3]}}\}$.
This should be followed by annihilation of the undesired one of the
two remaining irreps, using $[(3)]_6$.
Explicitly,
$$P_{[4,1,1]}=\frac{([(3)]_6+8)}{12}\cdot P^{\ast}$$
and
$$P_{[3,3]}=\frac{([(3)]_6-4)}{-12}\cdot P^{\ast}$$

The construction of the corresponding projection operators within $H_6(q)$
avoids the complication noted above, since all eigenvalues of the fundamental
invariant, $C_6$, are distinct.
Thus,
$$P_{[4,1,1]}(q)=\prod_{\Gamma\in{\cal{IRR}}_6\setminus \{[4,1,1]\} }
\frac{C_6-\Lambda_{C_6}^{\Gamma}}
{\Lambda_{C_6}^{[4,1,1]}-\Lambda_{C_6}^{\Gamma}}$$
and
$$P_{[3,3]}(q)=\prod_{\Gamma\in{\cal{IRR}}_6\setminus \{[3,3]\} }
\frac{C_6-\Lambda_{C_6}^{\Gamma}}
{\Lambda_{C_6}^{[3,3]}-\Lambda_{C_6}^{\Gamma}}$$

Taking the $q\rightarrow 1$ limit of $P_{[4,1,1]}(q)$ we observe that all the
factors except
$\frac{C_6-\Lambda_{C_6}^{[3,3]}}
{\Lambda_{C_6}^{[4,1,1]}-\Lambda_{C_6}^{[3,3]}}$
reduce to their $S_6$ counterparts.
To obtain the limiting form of this last factor we note that
$$\frac{C_6-\Lambda_{C_6}^{[3,3]}}{\Lambda_{C_6}^{[4,1,1]}-\Lambda_{C_6}^{[3,3]}}=
\frac{{\tilde{C}}_6-\Lambda_{{\tilde{C}}_6}^{[3,3]}}
{\Lambda_{{\tilde{C}}_6}^{[4,1,1]}-\Lambda_{{\tilde{C}}_6}^{[3,3]}}$$
Moreover,
$$\Lambda_{{\tilde{C}}_6}^{[4,1,1]} \sim 3\delta + 19\frac{\delta^2}{2!}
+\cdots$$
and
$$\Lambda_{{\tilde{C}}_6}^{[3,3]} \sim 3\delta + 7\frac{\delta^2}{2!} +\cdots$$
Therefore, using eq. \ref{equation:corres},
$${\tilde{C}}_6-\Lambda_{{\tilde{C}}_6}^{[3,3]} \sim \delta\Big( [(2)]_6-3\Big)
+\frac{\delta^2}{2!}\Big([(3)]_6+8\Big)+\cdots$$
and
$$\Lambda_{{\tilde{C}}_6}^{[4,1,1]}-\Lambda_{{\tilde{C}}_6}^{[3,3]}
\sim 12 \frac{\delta^2}{2!}+\cdots$$

The factor
$\frac{{\tilde{C}}_6-\Lambda_{{\tilde{C}}_6}^{[3,3]}}
{\Lambda_{{\tilde{C}}_6}^{[4,1,1]}-\Lambda_{{\tilde{C}}_6}^{[3,3]}}$
is multiplied by the projection operator
$$\prod_{\Gamma\in{\cal{IRR}}_6^{\ast} }
\frac{C_6-\Lambda_{C_6}^{\Gamma}}
{\Lambda_{C_6}^{[4,1,1]}-\Lambda_{C_6}^{\Gamma}}$$
that annihilates all irreps but ${\bf [4,1,1]}$ and ${\bf [3,3]}$, that
have a common eigenvalue of $[(2)]_6$.
Therefore, within that factor one can replace $[(2)]_6$ by its relevant
eigenvalue, 3,
to obtain
$$\frac{{\tilde{C}}_6-\Lambda_{{\tilde{C}}_6}^{[3,3]}}
{\Lambda_{{\tilde{C}}_6}^{[4,1,1]}-\Lambda_{{\tilde{C}}_6}^{[3,3]}} \sim
\frac{([(3)]_6+8)+O(\delta)}{12+O(\delta)}$$
or, finally,
$$P_{[4,1,1]} \equiv \lim_{q\rightarrow 1} P_{[4,1,1]}(q)= \frac{[(3)]_6+8}{12}
\prod_{\Gamma\in{\cal{IRR}}_6^{\ast} }
\frac{[(2)]_6-\lambda_{[(2)]_6}^{\Gamma}}
{\lambda_{[(2)]_6}^{[4,1,1]}-\lambda_{[(2)]_6}^{\Gamma}}$$
that coincides with the projection operator constructed directly within $S_6$,
utilizing both $[(2)]_6$ and $[(3)]_6$ to characterize the irreps of interest.

Clearly, if two or more irreps of $S_n$ were to have common eigenvalues for
$[(2)]_n$, $[(3)]_n$, $\cdots$, $[(k)]_n$, the relevant factor in the
projection operator
constructed within $H_n(q)$ in terms of $C_n$ would, in the limit $q\rightarrow
1$,
depend on $[(k+1)]_n$, the lowest single-cycle class-sum distinguishing between
these irreps.

\section{Traces of the Hecke algebra}

\hspace*{6 mm}
The problem of evaluating the traces of the Hecke algebra has received some
attention in recent years.
We refer to King and Wybourne \cite{King} for the  definition and the
discussion
of the properties of these traces, as well as for a presentation of the
pertinent earlier
references.

\noindent
By way of illustration of some properties that we shall need
let us demonstrate the observation made by King and Wybourne \cite{King},
according to
which the traces of elements of the same
connectivity class are equal.
In addition to the defining relations of the Hecke algebra, eq.
\ref{equation:Hecke},
we only use the property $tr(AB)=tr(BA)$.

\noindent
Thus, consider $tr(g_ig_{i+1}g_i)$, that we evaluate in two different ways
\begin{eqnarray}
tr(g_ig_{i+1}g_i) &=& tr(g_i^2g_{i+1})=(q-1)tr(g_ig_{i+1})+q\, tr(g_{i+1})
\nonumber \\
              &=&
tr(g_{i+1}g_ig_{i+1})=tr(g_ig_{i+1}^2)=(q-1)tr(g_ig_{i+1})+q\, tr(g_i)
\nonumber
\end{eqnarray}
Hence, $tr(g_i)=tr(g_{i+1})$.

\noindent
Similarly,
\begin{eqnarray}
tr(g_{i+2}g_ig_{i+1}g_{i+2}) &=& (q-1)tr(g_ig_{i+1}g_{i+2})+q\, tr(g_ig_{i+1})
\nonumber \\
                 &=& tr(g_ig_{i+2}g_{i+1}g_{i+2})=tr(g_ig_{i+1}g_{i+2}g_{i+1})=
tr(g_{i+1}g_ig_{i+1}g_{i+2})=\nonumber\\
                 & & {\phantom{spa}}
                     = tr(g_ig_{i+1}g_ig_{i+2})=tr(g_i^2g_{i+1}g_{i+2})=
\nonumber \\
& & {\phantom{spme space needed}} =(q-1)tr(g_ig_{i+1}g_{i+2})+q\,
tr(g_{i+1}g_{i+2})
\nonumber
\end{eqnarray}
from which it follows that $tr(g_ig_{i+1})=tr(g_{i+1}g_{i+2})$.
The generalization to Hecke algebra elements of any simply-connected class is
now
obvious.

As an example of a non-simply connected element consider
\begin{eqnarray}
tr(g_ig_{k+1}g_kg_{k+1})&=&tr(g_ig_kg_{k+1}g_k)=(q-1)tr(g_ig_kg_{k+1})+q\,
tr(g_ig_{k+1})
\nonumber \\
&=&(q-1)tr(g_ig_kg_{k+1})+q\, tr(g_ig_k) \nonumber
\end{eqnarray}
where $k>i+1$. This identity implies that $tr(g_ig_k)=tr(g_ig_{k+1})$, and
suggests
how similar connections can be established for arbitrary traces of terms
belonging
to common connectivity classes.

It is of some interest to note that the fundamental invariant determines these
traces as well.
To illustrate this point consider the simplest non-trivial case, corresponding
to
$H_3(q)$. In this case $C_3=g_1+g_2+\frac{1}{q}g_1g_2g_1$.
Obviously,
\begin{equation}
\label{equation:tr1}
tr(C_3)=3tr(g_1)+\frac{q-1}{q}tr(g_1g_2) \; .
\end{equation}
Furthermore,
$C_3^2=3qI+2(q-1)C_3+(\frac{1}{q}+1+q)(g_1g_2+g_2g_1+\frac{(q-1)}{q}g_1g_2g_1)$
so that
\begin{equation}
\label{equation:tr2}
tr(C_3^2)=3q \, tr(I)+(q^2+6q-6-\frac{1}{q})tr(g_1)
+(q^2+3q-2+\frac{3}{q}+\frac{1}{q^2})tr(g_1g_2)
 \; ,
\end{equation}
where $I$ is the unit operator.
Since both the eigenvalues of $C_3$ and the dimensionalities of the
corresponding irreps are known, equations \ref{equation:tr1}
and \ref{equation:tr2} can be solved for $tr(g_1)$ and $tr(g_1g_2)$.

Thus,
$$tr(g_1)=\frac{ (q^2+3q-2+\frac{3}{q}+\frac{1}{q^2})tr(C_3)
-\frac{q-1}{q}(tr(C_3^2)-3q \, tr(I))}
{2(q^2+2q+3+\frac{2}{q}+\frac{1}{q^2})}$$
and
$$tr(g_1g_2)=\frac{ -(q^2+6q-6-\frac{1}{q})tr(C_3)+3(tr(C_3^2)-3q \, tr(I))}
{2(q^2+2q+3+\frac{2}{q}+\frac{1}{q^2})}$$
For the one dimensional irrep ${\bf [3]}$ $\Lambda_{C_3}=q^2+2q$ so
$tr(I)=1$, $tr(C_3)=q^2+2q$, and $tr(C_3^2)=(q^2+2q)^2$.
Therefore, $tr(g_1)=q$ and $tr(g_1g_2)=q^2$.

\noindent
Similarly, for the two dimensional irrep ${\bf [2,1]} \;\; tr(I)=2$,
$tr(C_3)=2(q-1)$,
and \hfill\break
$tr(C_3^2)=2(q-1)^2$, so that $tr(g_1)=q-1$ and $tr(g_1g_2)=-q$,
and for the irrep ${\bf [1,1,1]} \;\; \hfill\break
tr(I)=1$, $tr(C_3)=-(2+\frac{1}{q})$, and
$tr(C_3^2)=(2+\frac{1}{q})^2$, hence $tr(g_1)=-1$ and $tr(g_1g_2)=1$, all in
agreement with King and Wybourne \cite{King}.

While this procedure can be followed to obtain the traces of higher Hecke
algebras,
a simpler procedure can be formulated in the following manner:

For a representation $\Gamma_n$ of $H_n(q)$ consider the basis states specified
by means of all the sequences of Young diagrams of the form
$\Gamma_2\Gamma_3\cdots\Gamma_n$, where each $\Gamma_{i+1}$ is obtained from
the
preceeding $\Gamma_i$ by the addition of one box. Since each such state is a
common eigenstate of the sequence of mutually commuting fundamental invariants
of
$H_2(q)\subset H_3(q)\subset \cdots\subset H_n(q)$, {\it i.e.,}
$C_2,\, C_3,\, \cdots,\, C_n$, it is also an eigenstate of the Murphy operators
\cite{Dipper} $L_2=C_2$, \hfill\break
$L_3=C_3-C_2,\, \cdots,\, L_n=C_n-C_{n-1}$.
Obviously, the eigenvalue of $L_i$, that we shall denote $\lambda_{L_i}$,
is equal to the q-content of the box added
to $\Gamma_{i-1}$ to form $\Gamma_i$. It is now a simple matter to construct
the (diagonal) representation matrices of the Murphy operators, from which
their
traces are readily obtained.

A straightforward computation yields
\begin{equation}
\label{equation:trC2}
tr(C_n)=\sum_{i=2}^n {n\choose i}\left(\frac{q-1}{q}\right)^{i-2}\tau_i
\end{equation}
where
$\tau_i=tr(g_1g_2\cdots g_{i-1})$.

Using eq. \ref{equation:trC2} we obtain
$$tr(L_i)=tr(C_i)-tr(C_{i-1})=\sum_{j=2}^i {{i-1}\choose{j-1}}
\left(\frac{q-1}{q}\right)^{j-2}\tau_j$$
or, inverting,
\begin{equation}
\label{equation:simply}
\tau_k=\left(\frac{q}{q-1}\right)^{k-2}
\sum_{i=0}^{k-2} (-1)^i{{k-1}\choose i} tr(L_{k-i})
\end{equation}

\noindent
As an illustration we note that in $H_3(q)$ we have:

\begin{enumerate}
\item
For the one-dimensional irrep ${\bf{\mbox{[3]}}}$ possessing the unique
sequence ${\mbox{\bf [1]}}\, {\mbox{\bf [2]}}\, {\mbox{\bf [3]}}$:
\hfill\break
$\lambda_{L_2}=q$,
$\lambda_{L_3}=q+q^2$.
Hence, using eq. \ref{equation:simply}, $tr(g_1)=q$ and $tr(g_1g_2)=q^2$.

\item
For the two-dimensional irrep ${\bf {\mbox{[2,1]}}}$, in the basis specified by
$$\left\{ \begin{array}{lllll}
{\mbox{\bf [1]}} & {\mbox{\bf [2]}} & {\mbox{\bf [2,1]}}:\;\;
& \lambda_{L_2}=q, & \lambda_{L_3}=-1 \\
{\mbox{\bf [1]}} & {\mbox{\bf [1,1]}} & {\mbox{\bf [2,1]}}:\;\;
& \lambda_{L_2}=-1, & \lambda_{L_3}=q
\end{array}\right.$$
$L_2$ and $L_3$ are represented by the
matrices
$$\left( \begin{array}{rr}
q & 0 \\
0 & -1
\end{array}\right)
\;\; {\mbox{and}} \;\;
\left( \begin{array}{rr}
-1 & 0 \\
0 & q
\end{array}\right), \; {\mbox{respectively.}}$$
Consequently, $tr(L_2)=tr(L_3)=q-1$ and, finally, $tr(g_1)=q-1$,
$tr(g_1g_2)=-q$.

\item
For ${\mbox{\bf [1]}}\, {\mbox{\bf [1,1]}}\, {\mbox{\bf [1,1,1]}}$:
$\lambda_{L_2}=-1$, $\lambda_{L_3}=-1-\frac{1}{q}$.
Hence, $tr(g_1)=-1$ and $tr(g_1g_2)=1$.
\end{enumerate}

\noindent
Similarly, in $H_4(q)$:

\begin{enumerate}
\item
For ${\mbox{\bf [1]}}\, {\mbox{\bf [2]}}\, {\mbox{\bf [3]}}\, {\mbox{\bf
[4]}}$:
$\lambda_{L_2}=q$,
$\lambda_{L_3}=q+q^2$ and $\lambda_{L_4}=q+q^2+q^3$.
Hence, $tr(g_1)=q$, $tr(g_1g_2)=q^2$, and $tr(g_1g_2g_3)=q^3$.

\item
For the three-dimensional irrep ${\bf {\mbox{[3,1]}}}$,
$$\left\{ \begin{array}{lllllll}
{\mbox{\bf [1]}} & {\mbox{\bf [2]}} & {\mbox{\bf [3]}} & {\mbox{\bf
[3,1]}}:\;\;
& \lambda_{L_2}=q, & \lambda_{L_3}=q+q^2, & \lambda_{L_4}=-1 \\
{\mbox{\bf [1]}} & {\mbox{\bf [2]}} & {\mbox{\bf [2,1]}} & {\mbox{\bf
[3,1]}}:\;\;
& \lambda_{L_2}=q, & \lambda_{L_3}=-1, & \lambda_{L_4}=q+q^2 \\
{\mbox{\bf [1]}} & {\mbox{\bf [1,1]}} & {\mbox{\bf [2,1]}} & {\mbox{\bf
[3,1]}}:\;\;
& \lambda_{L_2}=-1, & \lambda_{L_3}=q, & \lambda_{L_4}=q+q^2 \\
\end{array}\right.$$
Hence, $tr(L_2)=2q-1$, $tr(L_3)=q^2+2q-1$, and $tr(L_4)=2q^2+2q-1$, {\it,
i.e.,}
$tr(g_1)=2q-1$, $tr(g_1g_2)=q(q-1)$, and $tr(g_1g_2g_3)=-q^2$.
\end{enumerate}

\noindent
{\phantom{ele}}3.-5. etc.

To obtain the traces of the non-simply-connected element $g_1g_3$ we note that
\hfill\break
$L_2L_4=g_1(g_3+\frac{1}{q}g_2g_3g_2+\frac{1}{q^2}g_1g_2g_3g_2g_1)$,
so that \hfill\break
$tr(L_2L_4)=tr(g_1g_3)+\frac{q-1}{q}(q+\frac{1}{q})tr(g_1g_2g_3)+
2(q-1+\frac{1}{q})tr(g_1g_2)+(q-1)tr(g_1)$.

\noindent
Thus, for the irrep {\bf [4]} of $H_4(q)$, using the values of $\lambda_{L_2}$,
$\lambda_{L_4}$, $tr(g_1)$, $tr(g_1g_2)$ and $tr(g_1g_2g_3)$ obtained above,
$tr(g_1g_3)=q^2$.

\noindent
Similarly, for the irrep {\bf [3,1]}
$$L_2L_4\rightarrow \left(
\begin{array}{ccc}
-q & 0 & 0 \\
0 & q^2+q^3 & 0 \\
0 & 0 & -q-q^2
\end{array}
\right)$$
so that $tr(L_2L_4)=q^3-2q$ and $tr(g_1g_3)=q^2-2q$.

\noindent
etc.

Evaluating the trace of the product $L_2L_5$ we find that it can be expressed
in terms of $tr(g_1g_3)$, $tr(g_1g_3g_4)$, and traces of simply connected
terms.
Thus,
\begin{eqnarray}
tr(L_2L_5)&=&2\, tr(g_1g_3)+\left(\frac{q-1}{q}\right)tr(g_1g_3g_4) \nonumber
\\
&+&\left(\frac{q-1}{q}\right)^2\left(q+\frac{1}{q}\right)tr(g_1g_2g_3g_4)
+\left(\frac{q-1}{q}\right)\left(3q-2+\frac{3}{q}\right)tr(g_1g_2g_3)
\nonumber\\
&+&\left(3q-4+\frac{3}{q}\right)tr(g_1g_2)+(q-1)tr(g_1) \nonumber
\end{eqnarray}
Since $tr(g_1g_3)$, as well as all the simply-connected traces, have already
been
evaluated, $tr(L_2L_5)$ provides the next non-simply connected term,
$tr(g_1g_3g_4)$.

Continuing, we observe that products of pairs of Murphy operators provide the
traces
of all doubly-connected terms: $tr(L_2L_6)$ provides $tr(g_1g_3g_4g_5)$,
$tr(L_3L_6)$ provides $tr(g_1g_2g_4g_5)$, thus exhausting the doubly-connected
classes involving four Hecke generators, etc.
Triply-connected terms are generated by triple products of non-consecutive
Murphy operators. Thus, $tr(L_2L_4L_6)$ provides $tr(g_1g_3g_5)$, etc.
A more systematic study of the relations between traces of products of Murphy
operators
and traces of reduced non-simply connected terms will be presented elsewhere.

The following recursive procedure can be implemented to obtain the traces of
the
simply-connected
elements systematically:

For an irrep $\Gamma_n$ of $H_n(q)$ consider the set of irreps
$\Gamma_{n-1}\subset \Gamma_n$ of $H_{n-1}(q)$, obtained by eliminating
one box from $\Gamma_n$ in all possible ways.
Clearly,
$$tr(L_i)_{\Gamma_n}=\sum_{\Gamma_{n-1}\subset\Gamma_n}
tr(L_i)_{\Gamma_{n-1}}\;\;\;
i=2,\, 3,\, \cdots,\, n-1.$$
and
$$|\Gamma_n|=\sum_{\Gamma_{n-1}\subset\Gamma_n} |\Gamma_{n-1}|$$
where $|\Gamma_n|$ is the dimensionality of the irrep $\Gamma_n.$

\noindent
The trace of $L_n$ can now be evaluated very conveniently using either
$$tr(L_n)_{\Gamma_n}=
|\Gamma_n|\Lambda_{C_n}^{\Gamma_n} - \sum_{i=2}^{n-1} tr(L_i)_{\Gamma_n} \; ,$$
where $[\Gamma_n\setminus\Gamma_{n-1}]_q$ denotes the q-content of
the box that has been eliminated from $\Gamma_n$ to obtain $\Gamma_{n-1}$,
or the equivalent expression
$$tr(L_n)_{\Gamma_n}=
\sum_{\Gamma_{n-1}\subset\Gamma_n} |\Gamma_{n-1}|\,
[\Gamma_n\setminus\Gamma_{n-1}]_q\; .$$

This recursive procedure is ideally suited to implementation using symbolic
programming.
The traces of the Murphy operators of $H_2(q),\; H_3(q),\; \cdots,\; H_7(q)$,
evaluated (manually!) using this recursive procedure,
are presented in Tables 1-5.

Having obtained the traces of the Murphy operators, the traces of the
simply-connected
elements of the Hecke algebra can be evaluated using eq. \ref{equation:simply}.
For $H_2(q),\; H_3(q),\; \cdots,\; H_6(q)$ these traces agree with the values
presented by King and Wybourne \cite{King}. The traces of the simply-connected
terms in $H_7(q)$ are presented in Tables 6-7.

As illustrated above,
to obtain traces of non-simply-connected elements of the Hecke algebra we have
to
evaluate traces of products of Murphy operators.
For
$\prod_{i=1}^{\ell} L_{\alpha_i}$, \hfill\break
$2\leq\alpha_1\leq\alpha_2\leq\cdots\leq\alpha_{\ell}\leq n$,
we consider all the sequences of irreps leading to the desired irrep $\Gamma_n$
of
$H_n(q)$ from all possible irreps of $H_{\alpha_1-1}(q)$.
Let
$\Gamma_{\alpha_1-1}\subset\Gamma_{\alpha_1}\subset\Gamma_{\alpha_1+1}\subset\cdots
\subset\Gamma_n$ be such a sequence and let $\{\cdots\subset\Gamma_n\}$ denote
the
complete set of sequences.
Obviously,
$$tr\left(\prod_{i=1}^{\ell} L_{\alpha_i}
\right)=\sum_{\{\cdots\subset\Gamma_n\}}
|\Gamma_{\alpha_1-1}| \prod_{i=1}^{\ell}
[\Gamma_{\alpha_i}\setminus\Gamma_{\alpha_i-1}]_q \; .$$

\section{Correspondence between the eigenvalues of the fundamental invariants
of
$SU_q(N)$ and $H_n(q)$}

\hspace*{6 mm}
To introduce the relation between the fundamental invariant of the Hecke
algebra
$H_n(q)$ and that of
the quantum-unitary group $SU_q(N)$ we first recall that the space of an irrep
of
the latter is spanned by the Gelfand-Zetlin states \cite{Gelfand}

%\addtocounter{equation}{1}

\begin{picture}(425,140)
\put(0,0){\makebox(390,140){$
|h\rangle =
\left|
\begin{array}{ccccccccc}
h_{1,N} & & h_{2,N} &  & \cdots & & h_{N-1,N} & & h_{N,N} \\
 & h_{1,N-1} & & \cdots &  & \cdots & & h_{N-1,N-1} & \\
 & &\cdots & & \cdots & & \cdots & & \\
 & & & h_{1,2} & & h_{2,2} & & & \\
 & & & & h_{1,1} & & & &
\end{array}
\right. $}}
\put(395,70){\line(-1,3){18}}
\put(395,70){\line(-1,-3){18}}

%\put(415,0){\makebox(10,140){(14)}}

\end{picture}

For generic $q$ the parameters $h_{i,j}$ are integers satisfying the triangular
inequalities $h_{i,j+1}\geq h_{i,j} \geq h_{i+1,j+1} \;\;\; (i\leq j)$.
The top row specifies the irrep, $(h_{i,N}-h_{N,N})$ giving the length of the
$i$th
row of the Young diagram. A convenient convention is $h_{N,N}=0$.
The number of rows in a Young diagram corresponding to an irrep of $SU_q(N)$,
for generic $q$, is at most $N-1$, and the dimension is given by the number of
distinct
sets of integers $\{ h_{i,j} \; ; \; 1\leq i\leq j\leq N-1 \;\}$ satisfying the
triangular inequalities.

We will not reproduce the standard Drinfeld-Jimbo q-deformation of the $SU(N)$
algebra.
The matrix elements of the Chevalley triplets
$(h_k,\, e_k,\, f_k;\; k=1,\, 2,\, \cdots,\, N-1)$
are given by \cite{Chakrabarti,Arnaudon}
$$q^{\pm h_k}|h\rangle =
q^{\pm(2\sum_{i=1}^k h_{i,k} - \sum_{i=1}^{k+1} h_{i,k+1} - \sum_{i=1}^{k-1}
h_{i,k-1})}
|h\rangle $$
$$\langle
h_{j,k}-1|f_k|h\rangle=\langle h|e_k|h_{j,k}-1\rangle=
\left( -\frac{P_1(j,k) P_2(j,k)}{P_3(j,k)}\right)^{\frac{1}{2}}$$
where $|h_{j,k}-1\rangle$ differs from $|h\rangle$ only by the translation
$h_{j,k}\rightarrow (h_{j,k}-1)$, and
\begin{eqnarray}
P_1(j,k) &=& \prod_{i=1}^{k+1} [h_{i,k+1}-h_{j,k}-i+j+1]_s \nonumber \\
P_2(j,k) &=& \prod_{i=1}^{k-1} [h_{i,k-1}-h_{j,k}-i+j]_s \nonumber \\
P_3(j,k) &=& \prod_{1\leq i\leq k,\, i\neq j} [h_{i,k}-h_{j,k}-i+j+1]_s
[h_{i,k}-h_{j,k}-i+j]_s \nonumber
\end{eqnarray}
where $[x]_s$ is the symmetric q-bracket
$$[x]_s \equiv \frac{q^x-q^{-x}}{q-q^{-1}}$$

Define
$$ h_k= A_k^k-A_{k+1}^{k+1}$$
$$e_k=A_k^{k+1},\;\;\; f_k=A_{k+1}^k$$
$$A_{k-p}^{k+1}=A_{k-p}^{k+1-p}A_{k+1-p}^{k+1}-q^{-1}A_{k+1-p}^{k+1}A_{k-p}^{k+1-p}$$
$$A_{k+1}^{k-p}=A_{k+1}^{k+1-p}A_{k+1-p}^{k-p}-q^{-1}A_{k+1-p}^{k-p}A_{k+1}^{k+1-p}$$
The operator
$$ {\cal C}_2=\sum_{j>i} A_j^iA_i^j q^{(A_i^i+A_j^j-2j)}
+\frac{q}{(q-q^{-1})^2}\sum_{i=1}^N q^{2(A_i^i+i)} $$
commutes with all the generators.
Subtracting the invariant terms \cite{Chakrabarti}
\begin{equation}
\label{equation:eight}
\frac{1}{(q-q^{-1})^2}\left( q^{2\sum_{i=1}^N A_i^i+N(N+1)} + (N-1)\right)
\end{equation}
one obtains the correct classical $(q=1)$ limit
(this corresponds to a modified version of ${\cal{C}}_2$ of ref.
\cite{Chakrabarti}).
However, for our purposes it is more interesting to note the result
(obtained immediately by using the minimal state, that is annihilated by
$A_i^j$ for $j>i$) that
on the space of any irrep
\begin{equation}
\label{equation:nine}
(q-q^{-1})^2 q^{-(2N+3)}{{\cal{C}}_2} =
\left(\sum_{k=1}^N q^{2(h_{k,n}-k)}\right){\cal{I}}
\end{equation}
Defining $$ {{\hat{\cal{C}}}_2} \equiv
 (q-q^{-1})^2 q^{-(2N+3)}{{\cal{C}}_2} -q^{-2N}{\cal{I}}$$
we find that for $h_{N,N}=0$ the eigenvalue of ${\hat{\cal{C}}}_2$ on the irrep
$\Gamma$ is
$$ \Lambda_{{\hat{\cal{C}}}_2}^{\Gamma}
=\left( \sum_{k=1}^{N-1} q^{2(\ell_k -k)}\right){\cal{I}} $$
where $\ell_k$ is the length of the $k$th row of the Young diagram.
Note also that the linear invariant ${\cal{C}}_1 \equiv \sum_{i=1}^N A_i^i$ of
eq. \ref{equation:eight} has, for $h_{N,N}=0$, the eigenvalue
$$ \sum_{k=1}^{N-1} h_{k,N} =n $$
where $n=\sum_i \ell_i$ is the number of boxes in the Young diagram $\Gamma$
that labels the irrep.

Consider now the fundamental invariant of $H_n(q)$, introduced before.
The corresponding eigenvalues are ({\it cf.} eq. \ref{equation:eigen})
\begin{eqnarray}
\Lambda_{C_n}^{\Gamma} &=&
q\left(\sum_{(i,j)\in\Gamma}\frac{q^{j-i}-1}{q-1}\right) \nonumber \\
&=&\frac{q}{q-1}\left( \sum_{i\in\Gamma}\sum_{j=1}^{\ell_i}(q^{j-i}-1)\right)
\nonumber \\
&=&\frac{q}{q-1}\left(
\sum_{i\in\Gamma}\left(q^{-i+1}\frac{q^{\ell_i}-1}{q-1}-\ell_i\right)\right)
\nonumber
\end{eqnarray}
or,
$$ \left(\frac{q-1}{q}\right)^2 \Lambda_{C_n}^{\Gamma} =
\sum_{i\in\Gamma}\left( q^{\ell_i-i} -q^{-i} - \frac{q-1}{q} \ell_i\right) $$

If $\Gamma$ corresponds to the Young diagram that labels the irrep of $SU_q(N)$
considered
before than
\begin{eqnarray}
& & \sum_i \ell_i =n \nonumber \\
& & \sum_{i\in\Gamma} q^{-i} = \sum_{i=1}^{N-1} q^{-i}
= - \frac{q^{-N+1}-1}{q-1} \nonumber
\end{eqnarray}
Hence, if $q$ in $C_n(q)$ is replaced by $q^2$, than for irreps $\Gamma$ with
$n$ boxes
the precise correspondence between the eigenvalues of the
fundamental invariant of the
Hecke algebra and those of the quadratic Casimir of $SU_q(N)$ is
\begin{equation}
\label{equation:fifteen}
\left(\frac{q^2-1}{q^2}\right)^2 \Lambda_{C_n(q^2)}^{\Gamma}
+\frac{q^2-1}{q^2} n
= \Lambda_{{\hat{\cal{C}}}_2}^{\Gamma} +\frac{q^{2(-N+1)}-1}{q^2-1}
\end{equation}
Thus, through this relation $C_n(q^2)$ also characterizes the irreps of
$SU_q(N)$
corresponding to Young diagrams with $n$ boxes.

The relation between the transposition class-sum of the symmetric group and
the quadratic Casimir of $SU(N)$, that is the classical ($q=1$) limit of
eq. \ref{equation:fifteen}, was used by Gross \cite{Gross} to study the
$\frac{1}{N}$ expansion. For us the crucial terms are not the $n$ or $N$
dependent constants but
$$ \sum_k q^{2(\ell_k -k)} $$

The fact that one needs $C_n(q^2)$ in eq. \ref{equation:fifteen} can be viewed
in a more instructive way as follows. In the relation \ref{equation:Hecke}
$$g_i^2=(q-1)g_i+q$$
setting $q=(q^{\prime})^2$, $g_i=q^{\prime}g_i^{\prime}$ and then suppressing
the primes
for simplicity, one obtains
$$ g_i^2=(q-q^{-1})g_i +1 $$
This is the form used, for example, by Pan and Chen \cite{Pan}.
Using this form $C_n$ can be defined without explicit appearance of $q$
(or $q^{-1}$) and the correspondence with the $SU_q(N)$ Casimir is more direct.

For $SU_q(N)$ all the higher order Casimir operators are available
\cite{Reshetikhin,Bincer,Abdesselam}.
We will not reproduce them here. However, a few comments on the structure of
the eigenvalues can help to better appreciate the special role of
${\cal{C}}_2$ (or ${\hat{\cal{C}}}_2$).
The construction of Faddeev-Reshetikhin-Takhtadzhyan \cite{Reshetikhin}
(suitably normalized) can be shown to lead for the $p$th order invariant to the
eigenvalue
$$ \sum_{i_1<i_2<\cdots<i_p} q^{2\sum_{k=1}^p (h_{i_k,N}-i_k)} $$
This is not yet the most convenient form.
An operator giving the simplest symmetric polynomial of order $p$ in
$q^{h_{k,N}-k}$
namely,
$$ \sum_k q^{2p(h_{k,N}-k)} $$
would be desirable for relating to a corresponding invariant of $H_n(q)$.
But such operators remain to be constructed, both for $SU_q(N)$ and for
$H_n(q)$.
For $SU_q(N)$ the simplest, most compact construction of the Casimir operators
is that due to Bincer \cite{Bincer}, but the corresponding eigenvalues,
for all Casimirs but ${\cal{C}}_2$, are much more complicated.

Anyhow, since for generic $q\neq 1$, once $n$ is determined by the linear
invariant,
the fundamental Casimir ${\hat{\cal{C}}}_2$ completely specifies
the irreps, the necessity of the higher order Casimirs is not evident.

Since in eq. \ref{equation:nine}
$h_{k,N}\geq h_{k+1,N}$, so $ h_{k,n}-k >h_{k+1,N}-(k+1)\, . $
Hence, given an eigenvalue of ${\cal{C}}_2$,
simply arranging the powers of $q$ in decreasing order one gets an expression
of the type
$$ q^{2L_1}+q^{2L_2}+\cdots +q^{2L_{N-1}}\;\;\; (L_i>L_{i+1}) $$
Now,
$$ \ell_k=L_k+k\;\;\; (k=1,\, 2,\,\cdots,\, N-1) $$
completely specifies the Young diagram, giving the lengths of successive rows.
Again, since (setting $q^2=\exp(\delta)$)
$$ \sum_k q^{L_k} = 1+\delta(\sum_k L_k) +\frac{\delta^2}{2!}(\sum_k
L_k^2)+\cdots $$
the coefficient of $\delta^p$ gives the eigenvalue of a suitably defined $p$th
order classical invariant. Thus all the essential information is contained
in ${\cal{C}}_2$ (or ${\hat{\cal{C}}}_2$) in a simple way.
In fact, this is presumably true for the $q$-deformations of all algebras.

\section{Conclusions}

\hspace*{6 mm}
The fundamental invariant of the Hecke algebra $H_n(q)$ has been shown to fully
characterize the corresponding irreps.
In fact, when $q$ is neither zero nor a root of unity of order less than $n$,
the
eigenvalues of the fundamental invariant are polynomials in $q$ and
$\frac{1}{q}$
whose coefficients specify the structure of the corresponding Young diagram.
This is remarkable in view of the fact that the transposition class-sum, to
which this
fundamental invariant reduces in the limit $q\rightarrow 1$, exhibits (for
$n\geq 6$)
eigenvalue degeneracies.
The projection operators constructed in terms of the fundamental invariant
onto the various irreps are well behaved in the limit $q\rightarrow 1$, even
for
irreps that are not uniquely specified by the transposition class-sum.
These limiting projection operators generate the characters of all
the classes of $S_n$.

The eigenvalues of the quadratic Casimir operator of $SU_q(N)$ corresponding
to Young diagrams with $n$ boxes have been shown to be expressible in terms
of those of the fundamental invariant of $H_n(q)$. The former suffice to
specify
the irreps of $SU_q(N)$, {\it i.e.,} to determine the eigenvalues of all the
higher order
Casimir operators. We expect that the fundamental Casimir operators of other
quantum groups possess similar properties.

While the applications to the symmetric group and to the unitary and
quantum-unitary groups may be of mostly conceptual
interest, the results are demonstrated to lead to useful procedures for
evaluating
the traces of the Hecke algebra.

In this article we have considered only generic $q$. For $q$ a root of unity
the situation changes radically since the centre becomes suddenly much larger.
But then the formalism of ``fractional parts'' can be used
to unify the construction of the $SU_q(N)$ representations
\cite{Arnaudon,Abdesselam} (for $q$ generic and root of unity).
A similar approach to the $H_q(n)$ representations for $q$ a root of unity
might
be fruitful.
This aspect will be explored elsewhere.

\vspace{1cm}

\noindent
{\bf Acknowledgement}

We wish to acknowledge enlightening discussions with Daniel Arnaudon concerning
the
higher order $SU_q(N)$ invariants.

\end{document}